\documentclass[12pt]{amsart}
\usepackage{amsmath, amsthm, amscd, amsfonts}

\newtheorem{theorem}{Theorem}[section]
\newtheorem{lemma}[theorem]{Lemma}
\newtheorem{proposition}[theorem]{Proposition}
\theoremstyle{definition}
\newtheorem{definition}[theorem]{Definition}
\numberwithin{equation}{section}
\begin{document}
\title{$q$-norms are really norms}
\author{H. Belbachir, M. Mirzavaziri and M. S. Moslehian}
\address{Hac\`{e}ne Belbachir: USTHB/ Facult\'{e} de Math\'{e}matiques, B.P. 32, El Alia, 16111, Bab Ezzouar, Alger,
Alg\'{e}rie.} \email{hbelbachir@usthb.dz}
\address{Madjid Mirzavaziri: Department of Mathematics, Ferdowsi University, P. O. Box 1159, Mashhad 91775, Iran}
\email{mirzavaziri@math.um.ac.ir}
\address{Mohammad Sal Moslehian: Department of Mathematics, Ferdowsi University, P. O. Box 1159, Mashhad 91775, Iran}
\email{moslehian@ferdowsi.um.ac.ir} \subjclass[2000]{Primary
44B20; Secondary 46C05.} \keywords{Norm, q-norm, convexity}
\begin{abstract}
Replacing the triangle inequality, in the definition of a norm, by
$\left\| x + y \right\| ^{q}\leq 2^{q-1}\left( \left\| x\right\|
^{q} + \left\| y\right\| ^{q}\right) $, we introduce the notion
of a q-norm. We establish that every q-norm is a norm in the
usual sense,  and that the converse is true as well.
\end{abstract}
\maketitle

\section{Introduction}

The parallelogram law states that $\|x+y\|^2+\|x-y\|^2=2(\|x\|^2 +
\|y\|^2)$ holds for all vectors $x$ and $y$ in a Hilbert space.
This law implies that the so-called parallelogram inequality
$\|x+y\|^2\leq 2(\|x\|^2 + \|y\|^2)$ trivially holds. S. Saitoh
\cite{SAI} noted the inequality $\|x+y\|^2\leq
2(\|x\|^2+\|y\|^2)$ may be more suitable than the usual triangle
inequality. He used this inequality to the setting of a natural
sum Hilbert space for two arbitrary Hilbert spaces.

Obviously the classical triangle inequality in an arbitrary
normed space implies the above inequality. This motivates us to
introduce an apparently extension of the triangle inequality.
More precisely, we introduce the notion of a $q$-norm, by
replacing, in the definition of a norm, the triangle inequality
by $\| x+y\| ^{q}\leq 2^{q-1}\left( \| x\| ^{q}+\| y\|
^{q}\right)$, where $q \geq 1$. We establish that every q-norm is
a norm in the usual sense, and that the converse is true as well.
The reader is referred to \cite{J-L} for undefined terms and
notations.

\begin{definition} Let ${\mathcal X}$ be a real or complex linear space and $q \in [1, \infty)$. A mapping $%
\| \cdot \| :{\mathcal X}\rightarrow \left[ 0,\infty \right) $ is
called a $q$-norm on ${\mathcal X}$\ if it satisfies the following
conditions:

\begin{enumerate}
\item $\| x\| =0\Leftrightarrow x=0,$

\item $\| \lambda x\| =\| \lambda \|
\| x\| \ \ $for all $x\in {\mathcal X}$ and all scalar $\lambda ,$

\item $\| x+y\| ^{q}\leq 2^{q-1}\left( \|
x\| ^{q}+\| y\| ^{q}\right) \ $for all $x,y\in {\mathcal X}.$
\end{enumerate}
\end{definition}

We first prove a rather trivial result.

\begin{proposition}
Every norm in the usual sense is a $q$-norm.
\end{proposition}
\begin{proof} One can easily verify that the
function $f(t) = \frac{1 + t^q}{2} - (\frac{1 + t}{2})^q$ has a
nonnegative derivative and so it is monotonically increasing on
$[0, \infty)$. It follows that $(\frac{1 +
\frac{\|y\|}{\|x\|}}{2})^q \leq \frac{1 +
(\frac{\|y\|}{\|x\|})^q}{2}$ whenever $\|x\| \leq \|y\|$.
Therefore $\|\frac{x + y}{2}\|^q \leq (\frac{\|x\| + \|y\|}{2})^q
\leq \frac{\|x\|^q + \|y\|^q}{2}$ for all $x, y \in {\mathcal X}$.
It follows that $\|.\|$ is a $q$-norm.
\end{proof}

Now we state the following lemma which is interesting on its own
right.

\begin{lemma} Let ${\mathcal X}$ be a real or complex linear space. Let $%
\| \cdot \| :{\mathcal X}\rightarrow \left[ 0,\infty \right) $ be
a mapping satisfying (1) and (2) in the definition of a $q$-norm.
Then $\| \cdot \| $ is a norm if and only if the set $B=\left\{
x\mid \| x\| \leq 1\right\}$ is convex.
\end{lemma}

\begin{proof} If $\| \cdot \| $ is a norm, then $B$ is clearly a
convex set. Conversely, let $B$ be convex and $x,y\in {\mathcal
X}.$ We can assume that
$x \neq 0, y\neq 0$. Putting $x^{\prime }=\frac{x}{\| x\| }$ and $%
y^{\prime }=\frac{y}{\| y\| }$ we have $x^{\prime },y^{\prime
}\in B.$

Now $\lambda x^{\prime }+\left( 1-\lambda \right) y^{\prime }\in
B$ for all $0 \leq \lambda \leq 1.$ In particular, for $\lambda
=\frac{\|x\|}{\|x\| +\| y\| }$ we obtain%
\[
\| \frac{x}{\| x\| +\| y\| }+\frac{%
y}{\| x\| +\| y\| }\| =\| \lambda x^{\prime }+\left( 1-\lambda
\right) y^{\prime }\| \leq 1.
\]

So that $\| x+y\| \leq \| x\| +\| y\|.$
\end{proof}

We are just ready to prove our main result.

\begin{theorem}
Every $q$-norm is a norm in the usual sense.
\end{theorem}
\begin{proof} We shall show that $B = \left\{ x : \|
x\|
\leq 1\right\} $ is convex. Let $x,y\in B.$ Then we have%
\[
\| x+y\| ^{q}\leq 2^{q-1}\left( \| x\| ^{q}+\| y\| ^{q}\right)
\leq 2^{q-1}\left( 1+1\right) =2^{q}.
\]%
whence $\| \frac{x+y}{2}\| ^{2}\leq 1,$ so $\frac{1}{2}%
x+\left( 1-\frac{1}{2}\right) y\in B.$ Thus if $A=\left\{ \frac{k}{2^{n}}%
\mid n=1,2,\ldots ;k=0,1,\ldots ,n\right\}$, then for each
$\lambda \in A$ we have $\lambda x+\left( 1-\lambda \right) y\in
B.$

Let $0\leq \lambda \leq 1$ and $z=\lambda x+\left( 1-\lambda
\right) y.$ Since $A$\ is dense in $\left[ 0,1\right]$, there
exists a decreasing
sequence $\left\{ r_{n}\right\} $ in $A$ such that $\lim\limits_{n}r_{n}=%
\lambda .$ Put $\beta _{n}=\frac{1-r_{n}}{1-\lambda }.$ Obviously $0\leq
\beta _{n}\leq 1,$ $\lim\limits_{n}\beta _{n}=1$ and $\frac{r_{n}+\beta
_{n}-1}{r_{n}} \leq 1.$ Since $\frac{r_{n}+\beta _{n}-1}{r_{n}}x\in B$ and $%
r_{n}\in A$ we conclude that%
\[
\beta _{n}z = \lambda \beta _{n}x+\left( 1-\lambda \right) \beta _{n}y=r_{n}%
\frac{r_{n}+\beta _{n}-1}{r_{n}}x+\left( 1-r_{n}\right) y\in B.
\]%
Thus $\beta _{n} \| z\| =\| \beta _{n}z\| \leq 1$ for all $n.$
Tending $n$ to infinity we get $\| z\| \leq 1,$ i.e. $z\in B.$
\end{proof}

{\bf Acknowledgment.} We would like to sincerely thank Professor
Saburou Saitoh for his encouragement.

\end{document}